\documentclass{amsart}
\usepackage{amssymb}
\usepackage{amsmath}
\usepackage{amsfonts}

\theoremstyle{plain}

\numberwithin{equation}{section}
\input{tcilatex}

\begin{document}
\title[On Orthogonal Decomposition]{On Orthogonal Decomposition of $%
\tciLaplace ^{2}\left( \Omega \right) $}
\author{Dejenie A. Lakew}
\address{John Tyler Community College}
\email{dlakew@jtcc.edu}
\urladdr{http://www.jtcc.edu}
\date{February 26, 2015}

\begin{abstract}
\textit{In this short article we show an orthogonal decomposition of the
Hilbert space }$\tciLaplace ^{2}\left( \Omega \right) $ \textit{as }$%
\tciLaplace ^{2}\left( \Omega \right) =A^{2}\left( \Omega \right) \oplus 
\frac{d}{dx}\left( W_{0}^{1,2}\left( \Omega \right) \right) $, define
orthogonal projections and see some of its properties. We display some
decomposition of elementary functions as corollaries.
\end{abstract}

\maketitle

\textbf{Notations}.

\ \ \ 

\ \ \ \ 

Let $\Omega =\left[ 0,1\right] $

\ \ \ \ \ \ \ \ \ 

$\oplus $\ :\ \ \ \ \ Set direct sum

\ \ \ 

$\uplus $ :\ \ \ \ \ \ Unique direct sum of elements from mutually
orthogonal sets

\ 

$\left( \frac{d}{dx}\right) _{0}^{-2}$ : Inverse image of a second order
derivative of a traceless function

\ \ \ \ \ \ \ 

$A^{2}\left( \Omega \right) =\ker \frac{d}{dx}\cap \tciLaplace ^{2}\left(
\Omega \right) =\{f:\dint\limits_{\Omega }f^{2}dx<\infty \ni \left( \frac{d}{%
dx}\right) f=0$ on $\Omega \}$

$\Vert \ast \Vert :=\Vert \ast \Vert _{\tciLaplace ^{2}\left( \Omega \right)
}$

\ \ \ 

\ \ \ \ \ \ \ \ \ \ \ \underline{$\clubsuit $ \ \ \ \ \ \ \ \ \ \ \ \ \ \ \
\ \ \ \ \ \ \ \ \ \ \ \ \ \ \ \ \ \ \ \ \ \ \ \ \ \ \ \ \ \ \ \ \ \ \ \ \ \
\ \ \ \ \ \ \ \ }

\ \ \ \ \ \ \ 

We define the following function spaces

\ \ 

$(I)$ \ \ The Hilbert space of square integrable functions over $\Omega $ \ 
\begin{equation*}
\tciLaplace ^{2}\left( \Omega \right) =\{f:\Omega \longrightarrow 
\mathbb{R}
,\text{ measurable and \ }\dint\limits_{\Omega }f^{2}dx<\infty \}
\end{equation*}

\ 

$(II)$ \ The Sobolev space 
\begin{equation*}
W^{1,2}\left( \Omega \right) =\{f\in \tciLaplace ^{2}\left( \Omega \right)
:f_{w}^{\prime }\in \tciLaplace ^{2}\left( \Omega \right) \}
\end{equation*}%
where $f_{w}^{\prime }$ is a weak first order derivative of $f,i.e,$%
\begin{equation*}
\exists g\in \tciLaplace _{\text{loc}}\left( \Omega \right) :g=f_{w}^{\prime
}
\end{equation*}%
with 
\begin{equation*}
\dint\limits_{\Omega }g\varphi dx=-\dint\limits_{\Omega }f\varphi dx,\forall
\varphi \in C_{0}^{\infty }\left( \Omega \right)
\end{equation*}%
and

\ 

$(III)$ \ the traceless Sobolev space 
\begin{equation*}
W_{0}^{1,2}\left( \Omega \right) =\{f\in W^{1,2}\left( \Omega \right)
:f(0)=f(1)=0\}
\end{equation*}

The Hilbet space $\tciLaplace _{C^{1}}^{2}\left( \Omega \right) $ is an
inner product space with inner product 
\begin{equation*}
\langle ,\rangle _{\tciLaplace ^{2}\left( \Omega \right) }:\tciLaplace
^{2}\left( \Omega \right) \times \tciLaplace ^{2}\left( \Omega \right)
\longrightarrow 
\mathbb{R}%
\end{equation*}%
defined by 
\begin{equation*}
\langle f,g\rangle _{\tciLaplace ^{2}\left( \Omega \right)
}=\dint\limits_{\Omega }f\left( x\right) g\left( x\right) dx
\end{equation*}

and $W^{1,2}(\Omega )$ \ with an inner product 
\begin{equation*}
\langle f,g\rangle _{W^{1,2}(\Omega )}=\left( \langle f,g\rangle
_{\tciLaplace ^{2}\left( \Omega \right) }+\langle f_{w}^{\prime
},g_{w}^{\prime }\rangle _{\tciLaplace ^{2}\left( \Omega \right) }\right) ^{%
\frac{1}{2}}
\end{equation*}%
\ where $f_{w}^{\prime },g_{w}^{\prime }$ are weak first order derivatives.

\ 

\textbf{Definition 1}. \ \ For%
\begin{equation*}
f\in \tciLaplace ^{2}\left( \Omega \right) ,\text{ \ }\Vert f\Vert
_{\tciLaplace ^{2}\left( \Omega \right) }=\sqrt{\langle f,f\rangle
_{\tciLaplace ^{2}\left( \Omega \right) }}
\end{equation*}%
\ \ and for%
\begin{equation*}
f\in W^{1,2}\left( \Omega \right) ,\text{ \ }\Vert f\Vert _{W^{1,2}\left(
\Omega \right) }=\sqrt{\Vert f\Vert _{\tciLaplace ^{2}\left( \Omega \right)
}+\Vert f_{w}^{\prime }\Vert _{\tciLaplace ^{2}\left( \Omega \right) }}
\end{equation*}

\ \ \ \ \ \ \ \ \ \ \ \ \ 

With respect to the defined inner product above, we have the following
orthogonal decomposition

\ \ \ \ 

\ \ 

\textbf{Proposition 1}. \ (\textit{Orthogonal Decomposition)}\ \ \ \ 
\begin{equation*}
\ \tciLaplace ^{2}\left( \Omega \right) =A^{2}\left( \Omega \right) \oplus 
\frac{d}{dx}\left( W_{0}^{1,2}\left( \Omega \right) \right)
\end{equation*}

\ \ \ 

\ \ 

\textbf{Proof}. \ \ \ We need to show two things.

\ \ \ 

\ 

$(i)$ \ $A^{2}\left( \Omega \right) \oplus \frac{d}{dx}\left(
W_{0}^{1,2}\left( \Omega \right) \right) =\{0\}$

\ 

$(ii)$ $\ \forall f\in \tciLaplace ^{2}\left( \Omega \right) ,$ $\exists
g!\in A^{2}\left( \Omega \right) $ and $\exists h!\in \frac{d}{dx}\left(
W_{0}^{1,2}\left( \Omega \right) \right) $

\ \ \ \ \ \ 

such that%
\begin{equation*}
f=g\uplus h.
\end{equation*}

\bigskip

Indeed

\ 

$(i)$ \ Let $f\in A^{2}\left( \Omega \right) \cap \frac{d}{dx}\left(
W_{0}^{1,2}\left( \Omega \right) \right) $.

\ \ \ \ \ 

Then%
\begin{equation*}
f\in A^{2}\left( \Omega \right) \Longrightarrow \frac{d}{dx}f=0
\end{equation*}%
and so $f$ is a constant. Also%
\begin{equation*}
f\in \frac{d}{dx}\left( W_{0}^{1,2}\left( \Omega \right) \right)
\end{equation*}%
and hence 
\begin{equation*}
\exists h\in W_{0}^{1,2}\left( \Omega \right)
\end{equation*}%
such that%
\begin{equation*}
f=h^{\prime }
\end{equation*}%
But then as $f$ is a constant and we have 
\begin{equation*}
h=cx+d
\end{equation*}%
But%
\begin{equation*}
trh=0\text{ \ \ \ on \ }\partial \Omega =\{0,1\}
\end{equation*}%
and hence 
\begin{equation*}
h\left( 0\right) =0\Longrightarrow d=0
\end{equation*}%
and 
\begin{equation*}
h\left( 1\right) =0\Longrightarrow c=0
\end{equation*}%
Therefore 
\begin{equation*}
h\equiv 0\ \ \ \text{and hence \ \ }f\equiv 0.
\end{equation*}

\begin{equation*}
\therefore \text{ \ \ \ \ \ \ \ \ \ \ \ \ }A^{2}\left( \Omega \right) \cap 
\frac{d}{dx}\left( W_{0}^{1,2}\left( \Omega \right) \right) =\{0\}\ \ \ \ \
\ \ \ \ \ \ \ \ \ \ \ \ \ \ \ \ \ \ \ \ \ \ \left( \alpha \right)
\end{equation*}

$(ii)$ Let $f\in \tciLaplace ^{2}\left( \Omega \right) $. Then consider 
\begin{equation*}
\psi =\left( \frac{d}{dx}\right) _{0}^{-2}\left( \frac{d}{dx}\right) f
\end{equation*}%
\ which is in $W_{0}^{1,2}\left( \Omega \right) $ and let 
\begin{equation*}
g=f-\left( \frac{d}{dx}\right) \psi
\end{equation*}%
Then 
\begin{eqnarray*}
\frac{d}{dx}g &=&\frac{d}{dx}\left( f-\left( \frac{d}{dx}\right) \psi \right)
\\
&=&\frac{d}{dx}f-\frac{d^{2}}{dx^{2}}\left( \left( \frac{d}{dx}\right)
_{0}^{-2}\left( \frac{d}{dx}\right) f\right) \\
&=&0
\end{eqnarray*}%
Thus 
\begin{equation*}
g\in A^{2}\left( \Omega \right)
\end{equation*}%
and hence with 
\begin{equation*}
\eta =\left( \frac{d}{dx}\right) \psi \in \frac{d}{dx}\left(
W_{0}^{1,2}\left( \Omega \right) \right)
\end{equation*}%
we have 
\begin{equation*}
f=g\uplus \eta \ \ \ \ \ \ \ \ \ \ \ \ \ \ \ \ \ \ \ \ \ \ \ \ \ \ \ \ \ \ \
\ \ \ \ \left( \beta \right)
\end{equation*}

\ 

From $\left( \alpha \right) $ \ and $\left( \beta \right) $ \ follow the
proposition.

\ 

\textbf{Remark.} The suset $A^{2}\left( \Omega \right) $ is a closed set and
its orthogonal complement 
\begin{equation*}
\frac{d}{dx}\left( W_{0}^{1,2}\left( \Omega \right) \right) =\left(
A^{2}\left( \Omega \right) \right) ^{\bot }
\end{equation*}
as well.

\ \ \ \ \ \ \ \ 

Besides representation of elemnets of $\tciLaplace ^{2}\left( \Omega \right) 
$ is unique, i.e., 
\begin{equation*}
\forall f\in \tciLaplace ^{2}\left( \Omega \right) ,\exists !g\in
A^{2}\left( \Omega \right) \text{ \ and \ }!h\in \frac{d}{dx}\left(
W_{0}^{1,2}\left( \Omega \right) \right)
\end{equation*}
such that%
\begin{equation*}
f=g+h
\end{equation*}
which I denote it as 
\begin{equation*}
f=g\uplus h
\end{equation*}

\bigskip

\textbf{Definition 2.} \ Due to the orthogonal decomposition there are two
orthogonal projections 
\begin{equation*}
P:\tciLaplace ^{2}\left( \Omega \right) \longrightarrow A^{2}\left( \Omega
\right)
\end{equation*}%
\ and 
\begin{equation*}
Q:\tciLaplace ^{2}\left( \Omega \right) \longrightarrow \frac{d}{dx}\left(
W_{0}^{1,2}\left( \Omega \right) \right)
\end{equation*}%
with 
\begin{equation*}
Q=I-P
\end{equation*}%
where $I$ is the identity operator.

\ \ \ 

\ \ 

\textbf{Proposition 2}. \ \ $\ \ \forall f\in \tciLaplace ^{2}\left( \Omega
\right) $ \ \ we have \ \ 
\begin{equation*}
\left\langle P\left( f\right) ,Q\left( f\right) \right\rangle =0
\end{equation*}

\ \ \ 

\ \ 

\textbf{Proof.} \ \ Let $f\in \tciLaplace ^{2}\left( \Omega \right) $. Then%
\begin{equation*}
Pf\in A^{2}\left( \Omega \right)
\end{equation*}%
and so it is a constant and 
\begin{equation*}
Qf\in \frac{d}{dx}\left( W_{0}^{1,2}\left( \Omega \right) \right)
\end{equation*}%
and hence 
\begin{equation*}
\exists h\in W_{0}^{1,2}\left( \Omega \right)
\end{equation*}%
such that 
\begin{equation*}
Qf=h^{\prime }\ \ \ \ \text{with \ \ }trh=0
\end{equation*}%
Therefore 
\begin{equation*}
\left\langle P\left( f\right) ,Q\left( f\right) \right\rangle =\left\langle
P\left( f\right) ,h^{\prime }\right\rangle =\dint\limits_{\Omega
}P(f)h^{\prime }dx
\end{equation*}%
Then from integration by parts we have%
\begin{equation*}
\dint\limits_{\Omega }P(f)h^{\prime }dx=-\dint\limits_{\Omega }P(f)^{\prime
}hdx=0
\end{equation*}

\ 

since $P(f)$ is a constant and we have no boundary integral that might have
resulted from the application of integration by parts because of the
traceless of $h$.

\ 

\begin{equation*}
\therefore \ \ \text{\ \ \ \ \ \ \ \ \ \ }\left\langle P\left( f\right)
,Q\left( f\right) \right\rangle =0
\end{equation*}

\ 

\ 

\textbf{Proposition 3. }\ \ \ We have the following properties

\ \ \ 

\ 

$(i)$ \ \ \ $PQ=0$

\ 

$(ii)$ \ \ $P^{2}=P${}

\ 

$(iii)$ $\ Q^{2}=Q$ \ \ \ \ \ \ 

\ \ 

\ 

That is $P$ and $Q$ are \textit{idempotent}

\ \ 

\ \ 

\textbf{Proof. \ \ }Let $f\in \tciLaplace ^{2}\left( \Omega \right) $ and let%
\begin{equation*}
g=Pf\in A^{2}\left( \Omega \right)
\end{equation*}%
Then $g\in \tciLaplace ^{2}\left( \Omega \right) $ \ and let 
\begin{equation*}
\psi =\left( \frac{d}{dx}\right) _{0}^{-2}\left( \frac{d}{dx}g\right)
=\left( \frac{d}{dx}\right) _{0}^{-2}\left( 0\right)
\end{equation*}%
Then $\psi =0$ and setting%
\begin{equation*}
h=g-\underset{\underset{0}{\Vert }}{\underbrace{\frac{d}{dx}\psi }}
\end{equation*}%
\ we have 
\begin{equation*}
g=h+\underset{\underset{0}{\Vert }}{\underbrace{\frac{d}{dx}\psi }}
\end{equation*}%
\ with 
\begin{equation*}
Pg=h\text{ \ \ and \ }Qg=0
\end{equation*}

Therefore, 
\begin{equation*}
Pg=P^{2}f=h=g=Pf
\end{equation*}
and 
\begin{equation*}
Qg=QPf=0
\end{equation*}

Similarly let 
\begin{equation*}
\eta =Qf\in \frac{d}{dx}\left( W_{0}^{1,2}\left( \Omega \right) \right)
\end{equation*}

\ \ 

\textbf{Proof.} \ \ Let $f\in \tciLaplace ^{2}\left( \Omega \right) $. Then
we have the unique decomposition,

\begin{equation*}
f=Pf+Qf
\end{equation*}
But then 
\begin{eqnarray*}
\langle f,f\rangle &=&\langle Pf+Qf,Pf+Qf\rangle \\
&=&\langle Pf,Pf\rangle +\langle Qf,Qf\rangle
\end{eqnarray*}
That is 
\begin{equation*}
\Vert f\Vert ^{2}=\Vert Pf\Vert ^{2}+\Vert Qf\Vert ^{2}
\end{equation*}

\ \ 

We will look at few examples whose validity is supported from\textit{\
uniqueness} of representations in Hilbert spaces.

\ \ \ \ 

\ \ \ 

\textbf{Corollary 1.}

\ \ \ 

\ \ \ 

For $f\left( x\right) =x\in \tciLaplace ^{2}\left( \Omega \right) $ \ we have

\begin{equation*}
P\left( f\right) =\frac{1}{2}\ \text{\ \ and \ }Q\left( f\right) =x-\frac{1}{%
2}
\end{equation*}

and hence 
\begin{equation*}
f\left( x\right) =\frac{1}{2}\uplus \left( x-\frac{1}{2}\right)
\end{equation*}

\textbf{Proof. \ \ }Let 
\begin{equation*}
\psi =D_{0}^{-2}\left( Df\right) =\left( \frac{d}{dx}\right) _{0}^{-2}\left(
1\right) =\frac{1}{2}x^{2}-\frac{1}{2}x
\end{equation*}%
with 
\begin{equation*}
\frac{d}{dx}\psi =x-\frac{1}{2}
\end{equation*}%
and let 
\begin{equation*}
g=f-\frac{d}{dx}\psi =\frac{1}{2}
\end{equation*}%
Then 
\begin{equation*}
\frac{d}{dx}\left( g\right) =\frac{d}{dx}\left( f-\frac{d}{dx}\psi \right) =0
\end{equation*}%
and hence 
\begin{equation*}
f=g+\frac{d}{dx}\psi
\end{equation*}%
as a direct sum. That is%
\begin{equation*}
f=\frac{1}{2}\uplus \left( x-\frac{1}{2}\right)
\end{equation*}

\bigskip

\textbf{Corollary 2. \ \ \ }\ \ For \ $f\left( x\right) =x$ $\ \ \ $%
\begin{equation*}
\left\langle P\left( f\right) ,Q\left( f\right) \right\rangle =0
\end{equation*}

\ \ 

\ \ \ 

\textbf{Proof}. \ \ \ Indeed 
\begin{eqnarray*}
\ \left\langle P\left( f\right) ,Q\left( f\right) \right\rangle
&=&\dint\limits_{\Omega }\frac{1}{2}\left( x-\frac{1}{2}\right) dx \\
&=&\frac{1}{2}\left( \frac{x^{2}}{2}-\frac{x}{2}\right) _{0}^{1} \\
&=&0
\end{eqnarray*}

\ 

\textbf{Corollary 3.} \ \ $\Vert x\Vert ^{2}=\Vert \frac{1}{2}\Vert
^{2}+\Vert x-\frac{1}{2}\Vert ^{2}$

\ \ \ \ \ 

\textbf{Corollary 4. \ \ }For $f\left( x\right) =x^{2}$%
\begin{equation*}
P\left( f\right) =\frac{1}{3}\ \ \ \ \text{and \ }\ Q\left( f\right) =x^{2}-%
\frac{1}{3}
\end{equation*}

\ \ \ 

\ \ 

\textbf{Proof}. \ \ \ Let 
\begin{equation*}
\psi =\left( \frac{d}{dx}\right) _{0}^{-2}\left( \frac{d}{dx}f\right)
=\left( \frac{d}{dx}\right) _{0}^{-2}\left( 2x\right)
\end{equation*}%
\begin{equation*}
\Longrightarrow \text{ \ \ \ \ \ \ \ }\psi \left( x\right) =\frac{1}{3}x^{3}-%
\frac{1}{3}x
\end{equation*}%
and let 
\begin{eqnarray*}
g &=&f-\frac{d}{dx}\psi \\
&=&x^{2}-\left( x^{2}-\frac{1}{3}\right) \\
&=&\frac{1}{3}
\end{eqnarray*}%
and so 
\begin{equation*}
g\in \ker \frac{d}{dx}
\end{equation*}%
\ and so 
\begin{equation*}
f=g\uplus \frac{d}{dx}\psi =\frac{1}{3}\uplus \left( x^{2}-\frac{1}{3}\right)
\end{equation*}%
which signifies 
\begin{equation*}
P\left( f\right) =\frac{1}{3}\ \ \ \ \text{and \ }Q\left( f\right) =x^{2}-%
\frac{1}{3}
\end{equation*}%
\ with 
\begin{equation*}
\left\langle \frac{1}{3},\text{ }\left( x^{2}-\frac{1}{3}\right)
\right\rangle =0
\end{equation*}

\ \ \ 

\ 

\textbf{Corollary 5.} \ \ \ \ $\Vert x^{2}\Vert ^{2}=\Vert \frac{1}{3}\Vert
^{2}+\Vert \left( x^{2}-\frac{1}{3}\right) \Vert ^{2}$

\ 

\textbf{Proposition 4}. \ \ For the orthogonal projections $P$ and $Q$ we
have the following results

\ \ \ \ \ \ 

\ \ 

$(i)$ $\ \ \ x^{n}=\frac{1}{n+1}\uplus \left( x^{n}-\frac{1}{n+1}\right) $

\ \ \ 

i.e. 
\begin{equation*}
P(x^{n})=\frac{1}{n+1},\text{ \ }Q\left( x^{n}\right) =x^{n}-\frac{1}{n+1}
\end{equation*}

\ 

$(ii)$ $\ \ \ e^{x}=\left( e-1\right) \uplus \left( e^{x}+1-e\right) $

\ 

i.e., 
\begin{equation*}
P\left( e^{x}\right) =e-1,\text{ \ \ }Q\left( e^{x}\right) =e^{x}+1-e
\end{equation*}

\ 

$(iii)$ $\ \ P(\cos x)=\sin 1$, $\ Q\left( \cos x\right) =\cos x-\sin 1$

\ 

so that 
\begin{equation*}
\cos x=\sin 1\uplus \left( \cos x-\sin 1\right)
\end{equation*}

\ 

$(iv)$ $\ \ P(\sin x)=1-\cos 1$, $\ Q\left( \sin x\right) =\sin x+\cos 1-1$
\ 

\ \ 

so that 
\begin{equation*}
\sin x=\left( 1-\cos 1\right) \uplus \left( \sin x+\cos 1-1\right)
\end{equation*}

\textbf{Proof of \ }$\mathbf{(iii)}$. \ \ Let 
\begin{equation*}
\psi =\left( \frac{d}{dx}\right) _{0}^{-2}\left( \frac{d}{dx}\cos x\right)
=\sin x-\left( \sin 1\right) x
\end{equation*}%
\begin{equation*}
\Longrightarrow \text{ \ \ \ \ \ }\frac{d}{dx}\psi \left( x\right) =\cos
x-\sin 1
\end{equation*}%
Then set 
\begin{equation*}
g=f-\frac{d}{dx}\psi =\sin 1\in \ker \frac{d}{dx}
\end{equation*}

Thus 
\begin{equation*}
\cos x=\sin 1\uplus \left( \cos x-\sin 1\right)
\end{equation*}

and hence 
\begin{equation*}
P\left( \cos x\right) =\sin 1\ \ \ \ \text{and \ \ }Q\left( \cos x\right)
=\cos x-\sin 1
\end{equation*}

\ 

\textbf{Corollary 6}. \ \ 

\ \ \ \ 

\ \ 

$(i)$ $\ \ \Vert x^{n}\Vert ^{2}=\Vert \frac{1}{n+1}\Vert ^{2}+\Vert x^{n}-%
\frac{1}{n+1}\Vert ^{2}$

\ 

$(ii)$ $\ \Vert e^{x}\Vert ^{2}=\Vert e-1\Vert ^{2}+\Vert e^{x}+1-e\Vert
^{2} $

\ 

$(iii)$ $\Vert \cos x\Vert ^{2}=\Vert \sin 1\Vert ^{2}+\Vert \cos x-\sin
1\Vert ^{2}$

\ \ 

\ \ \ \ 

$\clubsuit $\underline{ \ \ \ \ \ \ \ \ \ \ \ \ \ \ \ \ \ \ \ \ \ \ \ \ \ \
\ \ \ \ \ \ \ \ \ \ \ \ \ \ \ \ \ \ \ \ \ \ \ \ \ \ \ }

\ \ \ 

\ \ 

\textbf{References}.

\ 

$\left[ 1\right] $. \ Clifford analytic complete function systems for
unbounded domain, Math. Meth. Appl. Sci., Vol. 25 (2002) 1527-1539, Dejenie
A. Lakew and John Ryan.

$\left[ 2\right] $. \ Complete function systems and decomposition results
arising in Clifford analysis, Comp. Meth. Fun. Theory, CMFT 2, No. 1 (2002)
215-228, Dejenie A. Lakew and John Ryan.

$\left[ 3\right] $. \ Clifford analysis over Orlicz-Sobolev spaces,
arXiv:1409.8380v1, Dejenie A. Lakew and Mulugeta Alemayehu.

$\left[ 4\right] .$ On a direct decomposition of the space $\tciLaplace
^{p}\left( \Omega \right) $, Z. Anal. Anwend. 18 (1999) 839-884, U. K\"{a}%
hler.

\end{document}